\documentclass[red,11pt,a4paper]{article}
\usepackage{fullpage}

\usepackage{cite}

\usepackage{amsmath}
\usepackage{amscd}
\usepackage{amssymb}
\usepackage{bm}
\usepackage{graphicx}
\usepackage{latexsym}
\usepackage{color}
\usepackage{verbatim}

\bf
\newcommand{\vr}{\varrho}
\newcommand{\vrd}{\vr_\delta}
\newcommand{\vrs}{\varrho^*}
\newcommand{\pt}{\partial_{t}}

\newcommand{\Ov}[1]{\overline{#1}}
\newcommand{\lap}{\Delta}

\newcommand{\Div}{{\rm div}}
\newcommand{\Dt}{\frac{ d}{dt}}
\newcommand{\vu}{\vc{u}}
\newcommand{\vud}{\vu_\delta}
\newcommand{\vc}[1]{{\bf #1}}
\newcommand{\vS}{{\bf S}}
\newcommand{\Grad}{\nabla}
\newcommand{\R}{\mathbb{R}}
\newcommand{\lr}[1]{\left( #1 \right)}
\newcommand{\ep}{\varepsilon}
\newcommand{\QTb}{[0,T]\times \overline{\Omega}}
\newcommand{\eq}[1]{\begin{equation}
\begin{split}
#1
\end{split}
\end{equation}}
\newcommand{\eqh}[1]{\begin{equation*}
\begin{split}
#1
\end{split}
\end{equation*}}
\newcommand{\bo}{| _{\partial\Omega}}
\newcommand{\intO}[1]{\int_{\Omega} #1 \ \dx}

\newcommand{\intT}[1]{\int_0^T #1 \ \dt}
\newcommand{\intTO}[1]{\int_0^T\!\!\!\! \int_{\Omega} #1 \ \dxdt}
\newcommand{\intTOB}[1]{ \int_0^T\!\!\!\! \int_{\Omega} \left( #1 \right) \ \dxdt}
\newcommand{\dx}{{\rm d} {x}}
\newcommand{\dt}{{\rm d} t }
\newcommand{\dxdt}{\dx \ \dt}
\newtheorem{thm}{Theorem}
\newtheorem{lemma}[thm]{Lemma}

\newtheorem{df}{Definition}
\newtheorem{rmk}{Remark}

\title{\bf Transport of congestion in two-phase compressible/incompressible flows}
\author{Pierre Degond$^1$, Piotr Minakowski$^2$, Ewelina Zatorska$^{1}$}

\begin{document}
 \maketitle 
\centerline{1. Department of Mathematics, Imperial College London, }
\centerline{London SW7 2AZ, United Kingdom.}
\bigskip
  \centerline{2.  Interdisciplinary Center for Scientific Computing, Heidelberg University,}
 \centerline{Im Neunheimer Feld 205, D-69120 Heidelberg, Germany.}
 \bigskip
 
\noindent {\small {\bf Abstract:} We study the existence of weak solutions to the two-phase model of crowd motion. The model encompasses the flow in the uncongested regime (compressible) and the congested one (incompressible) with the free boundary separating the two phases. The congested regime appears when the density in the uncongested regime $\vr(t,x)$ achieves a threshold value $\vrs(t,x)$ that describes the comfort zone of individuals. This quantity is prescribed initially and transported along with the flow. We prove that this system can be approximated by the fully compressible Navier-Stokes system with a singular pressure, supplemented with transport equation for the congestion density. We also present the application of this approximation for the purposes of numerical simulations in the one-dimensional domain.}

 \bigskip
 
 
 \section{Introduction} 

Mathematical modelling of crowd dynamics is a very challenging problem. The recently rapidly growing literature on this subject covers several parallel approaches. We can distinguish, for example, the mean-field game models \cite{Do10, LW11}, in which the individuals behave as the players following some strategy, or optimizing certain cost; the microscopic models which describe precise position and velocity of an individual (Individual-Based-Models) using Newtonian framework \cite{GH08, HM97,TCP06,HFV00}; or the macroscopic models formulated in the language developed for the fluids \cite{CMW16,He92,DeFu14,Hu03,PiTo11}. The behaviour of the crowd in the later is characterised by some averaged quantities such as the number density or mean velocity. The macroscopic models, although less precise than the microscopic ones, are computationally more affordable. Moreover, they allow for asymptotic studies that proved to be useful for understanding various aspects like: swarming or pattern formation observed in the experiments. 

Our aim is to analyze the  free-boundary two-phase fluid system that could be used to model the congestions in the large group of individuals in a bounded area. Our individuals do not follow any neighbour trying to align their velocity or reach a certain evacuation point. They are just the agents that have their individual preferences for how close they let the closest neighbour to approach and they carry this information with them in the course of motion. We prescribe their initial velocity that determines their desired direction of motion and check how the individual preferences as well as the initial distribution of the agents determines creation of congestions.

Our system writes as follows:
 \begin{subequations}\label{sysSl}
 \begin{equation}\label{rho}
 \pt\vr+ \Div (\vr\vu) = 0, 
 \end{equation}
\begin{equation}\label{mom}
 \partial_t (\vr\vu) + \Div (\vr\vu \otimes \vu) + \Grad\pi +\Grad p\lr{\frac{\vr}{\vr^*}}-\Div\vc{S}(\vu) = \vc{0},
 \end{equation}
\begin{equation}\label{rho_star}
 \pt\vr^*+\vu\cdot\Grad\vr^*=0,
 \end{equation}
 \begin{equation}\label{cons0}
 0\leq \vr\leq\vrs,
 \end{equation}
\begin{equation}\label{div0}
  \Div\vu =0 \ \text{in} \ \{\vr=\vrs\},
 \end{equation}
\begin{equation}
  \pi\geq 0\ \mbox{in} \ \{\vr=\vrs\},\quad \pi= 0\ \mbox{in} \ \{\vr<\vrs\}.\label{pineq0}
  \end{equation}
 \end{subequations}
with the unknowns: $\vr=\vr(t,x)$ -- the mass density, $\vu=\vu(t,x)$ -- the velocity vector field, $\vrs=\vrs(t,x)$ -- the congestion density, also referred to as the barrier or the threshold density, and $\pi$ -- the congestion pressure appearing only when $\vr=\vrs$.
 
The barotropic pressure is an explicit function of $\frac{\vr}{\vr^*}$
\eq{\label{pbar}
p\left(\frac{\vr}{\vr^*}\right)=\lr{\frac{\vr}{\vr^*}}^\gamma,\quad \gamma>1,
}
and plays the role of the background pressure.

The stress tensor $\vS$ is a known function of $\vu$, characteristic for the Newtonian fluid, namely
\eq{\label{S}
\vS=\vS(\vu)= 2\mu\, {\rm D}(\vu)+\lambda\Div\vu\, {\rm I},\quad \mu>0, \ 2\mu+\lambda>0,
}
where ${\rm D}(\vu)=(\Grad\vu+\Grad^T\vu)/2$ denotes the symmetric part of the gradient of $\vu$, and ${\rm I}={\rm I}_3$ is the identity matrix. 

It is justified to call the above system the two-phase system because for $\vr(t,x)<\vrs(t,x)$ it behaves as the compressible Navier-Stokeas system with the barotropic pressure, while when the congestion is achieved, i.e. for $\vr(t,x)=\vrs(t,x),$ the system behaves like the incompressible Navier-Stokes equations. We thus observe a switching between two phases: compressible and incompressible depending on the size of the density ratio $\frac{\vr}{\vrs}$. 
The fluid systems with congestion constraints have been recently intensively studied, especially in the hyperbolic regime \cite{BoBrCoRi, Be02, Berthelin16}. The first analytical result for system \eqref{sysSl} with $\vr^*=1$ is due to P.-L. Lions and N. Masmoudi \cite{LM99}, who showed that it can be obtained as a limit of compressible Navier-Stokes equations with barotropic pressure $\vr^{\gamma}$ with $\gamma\to \infty$, similar studies were performed recently for the model of tumour growth \cite{PeVa15}.

In the system \eqref{sysSl} $\vr^*$ models preferences of the individuals, it is given initially and then transported with the flow. Therefore, $\vr^*$ depends on time and position, but more importantly it depends on initial configuration $\vrs_0$. The form of $\vrs$ relaxes the restrictions from the models studied in \cite{BrPeZa, PeZa}, where the threshold density $\vrs$ was either assumed to be constant or independent of time. This allows to cover more physical applications. Including the transport of the congestion density $\vrs$ allows also to study the system \eqref{sysS} with the contribution from the pressure in the form of the pure gradient, without factor $\vrs$ as it was done in \cite{PeZa}.

We will consider the system \eqref{sysSl} in the 3-dimensional domain $\Omega$ with the smooth boundary $\partial \Omega$, and the Dirichlet boundary conditions for the velocity vector field
\eq{\label{boundary}
\vu\bo=\vc{0}.}

The initial conditions are given by:
\eqh{\vr(0,x)=\vr_0(x),\quad(\vr\vu)(0,x)=\vc{m}_0(x),\quad \vrs(0,x)=\vrs_0(x),}
and we assume that they satisfy:
\begin{itemize}
\item $\vr_0\geq 0$, $\intO{\vr_0}>0,$
\item $\vc{m}_0 =\bf{0}\  \mbox{a.e. in } \Omega$, $ \vc{m}_0\in L^2(\Omega)$,
\item $\vr_0\leq\vrs_0,  \  \mbox{a.e. in }\Omega$, $\vr_0\neq\vrs_0$, $\vrs_0\in L^\infty(\Omega)$.
\end{itemize}
 Moreover, we assume that in the region of the absence of the individuals $\vr_0(x)=0$, the congestion density is equal to a constant value, being the characteristic mean preference of the group:
 \eqh{\vrs_0\big|_{\{\vr_0=0\}}=\tilde{\vrs}>0.}
 
 The main result of this paper is the existence of solutions to the system \eqref{sysSl} under the aforementioned assumptions on the constitutive relations and the initial condition, in the sense of the following definition.
 
\begin{df}[Weak solution]\label{Def1}
A quadruple $(\vr,\vu,\vrs,\pi)$ is called a weak solution to \eqref{sysSl}  if equations \eqref{rho}, \eqref{mom}, \eqref{rho_star}
are satisfied in the sense of distributions, the constraints \eqref{cons0} satisfied a.e. in $(0,T)\times \Omega$, the divergence free condition \eqref{div0} is satisfied a.e. in $\{\vr=\vrs\}$, 
and the following regularity properties hold
\begin{align*}
&\vr\in C_w([0,T];L^\infty(\Omega)), \\
&\vrs\in C_w([0,T];L^\infty(\Omega)),\\
& \vu \in L^2(0,T;W^{1,2}_0(\Omega, \R^3)),\quad \vr|\vu|^2 \in L^{\infty}(0,T; L^1(\Omega)),\\
&\pi\in {\cal M}^+ ((0,T)\times \Omega).
\end{align*}
Moreover, $\pi$ is sufficiently regular so that the condition
\eqh{(\vrs-\vr)\pi=0,}
is satisfied in the sense of distributions on $(0,T)\times \Omega$.
\end{df}

The main theorem of the paper states as follows.
 \begin{thm}\label{Th:main}
Let the initial conditions $\vr_0, \ \vc{m_0}, \ \vrs_0$ satisfy the conditions above. Then the system \eqref{sysSl} with $p$ and $\vS$ given by \eqref{pbar}, \eqref{S} respectively, has a weak solution in the sense of Definition \ref{Def1}.
\end{thm}
 
\begin{rmk}
 The same result holds for the lower dimensions, d=1,2. 
 \end{rmk}
The core  of the proof of Theorem \ref{Th:main} is to show that the system \eqref{sysSl} can be obtained as a limit when $\ep\to 0$ of the following approximation
\begin{subequations} \label{sysS}
 \eq{\label{sysSa}
  \pt\vr+ \Div (\vr\vu) = 0}
 \eq{\label{sysSb}
 \partial_t (\vr\vu) + \Div (\vr\vu \otimes \vu) + \Grad\pi_\ep\lr{\frac{\vr}{\vr^*}}+\Grad p\lr{\frac{\vr}{\vr^*}} -\Div\vc{S}(\vu) = \vc{0}}
 \eq{\label{sysSc}
 \pt\vr^*+\vu\cdot\Grad\vr^*=0,}
\end{subequations}
where the $\pi_\ep$ stands for the  singular pressure of the form
\eq{\label{pie}
\pi_\ep(r)=\ep \frac{\lr{\frac{\vr}{\vr^*}}^\alpha}{\lr{1-\frac{\vr}{\vr^*}}^\beta}, \quad \alpha\geq0,\ \beta>0.
}
A similar form of the pressure
\[\displaystyle\ep\nabla \frac{1}{\left(\frac{1}{\vr}-\frac{1}{\vr^*}\right)^\beta}\] 
was proposed in {\rm \cite{BeBr}} or {\rm \cite{DeHuNa}}. Singularities of the pressure of this type were also previously studied in the context of traffic models \cite{BeBr, Berthelin_M3AS_08, Berthelin_ARMA_08}, collective dynamics \cite{DeHuNa, DeHu}, or granular flow \cite{Maury07, Perrin16}.

Our first goal is to prove the existence of solutions to a certain reformulation of the system \eqref{sysS} for $\ep>0$ fixed. Similar system, but with the barotropic form of the pressure $\pi=\vr^\gamma$ was studied recently in \cite{M3NPZ}, see also the stability result from \cite{Mich15} and the low Mach number analysis from \cite{FeKlNoZa}. It was proved that for certain values of parameter $\gamma$, there is an equivalence between three different definitions of weak solutions to \eqref{sysS} based on the reformulations of the transport equation for the entropic variable. 

Here we will essentially use the formulation involving a new unknown $Z$, that for $\vr,\ \vrs$ smooth enough can be identified with the density fraction $Z=\frac{\vr}{\vrs}$. We can formally check, that dividing  \eqref{sysSa} by $\vrs$, multiplying \eqref{sysSc} by $-\frac{\vr}{{\vrs}^2}$, and summing the resulting expressions, the system \eqref{sysS} can be transformed to the following one
\begin{subequations} \label{sysSZ}
 \eq{\label{sysSZa}
 \pt\vr+ \Div (\vr\vu) = 0}
  \eq{\label{sysSZb}\partial_t (\vr\vu) + \Div (\vr\vu \otimes \vu) + \Grad\pi_\ep(Z) +\Grad p(Z)-\Div\vc{S}( \vu) = \vc{0}}
  \eq{\label{sysSZc}\pt Z+\Div(Z\vu)=0}
  \end{subequations}
with the initial data
\eq{\label{initialZ}
\vr(0,x)=\vr_0(x),\quad(\vr\vu)(0,x)=\vc{m}_0(x),\quad Z(0,x)=Z_0(x).}
In consistency with the assumptions on $\vr_0,\ \vc{m}_0,\ \vrs_0$ from the previous section, we postulate that $\vr_0,\ \vc{m}_0,\ Z_0$ satisfy
\begin{equation}\label{initialdata2}
\begin{array}{c}
\displaystyle 0 \leq c_\star \vr_0 \leq Z_0 \leq c^\star \vr_0 \mbox{ a.e. in } \Omega,  \quad \mbox{for}\quad 0<c_\star\leq c^\star <\infty,\\
 0<\intO{\vr_0},\quad Z_0\leq1,\quad  \intO{Z_0}< |\Omega|,\quad \frac{\vr_0}{Z_{0}}\big|_{\{\vr_0=0\}}=\tilde{\vrs},\\
 \displaystyle \vc{m}_0 \in L^2(\Omega,\R^3),\quad \vc{m}_0 =0\  \mbox{a.e. in } \{\vr_0=0\}.\\
\end{array}
\end{equation}
The first condition from \eqref{initialdata2} should be understood as the restriction of the initial congestion density, having in mind our notation $Z=\frac{\vr}{\vr^*}$ we see that the condition above means that $\vrs_0$ cannot be zero on the regions the density $\vr_0$ in positive. The left condition means that the congestion density is initially bounded. The restriction on the integral $\intO{Z}<|\Omega|$ means that the assumption $\vr_0\neq\vrs_0$ holds on a set of non-zero Lebesgue measure.

\begin{df}[Weak solution of the approximate system]\label{weaksolZ}
Suppose that the initial conditions satisfy (\ref{initialdata2}).
We say that the triplet  $(\vr,\vu,Z)$ is a weak solution of problem \eqref{sysSZ} with the initial and boundary conditions (\ref{initialZ}), (\ref{boundary}) if
\begin{equation*}
(\vr,\vu,Z) \in  L^\infty((0,T)\times \Omega) \times L^2(0,T;W^{1,2}_0(\Omega,\R^3)) \times L^\infty((0,T)\times\Omega),
\end{equation*}
and for any $T>0$ we have:
\begin{description}
\item{(i)}
$ \vr \in C_w([0,T];L^\infty(\Omega))$, and \eqref{sysSZa} is satisfied in the weak sense
\eq{
\intO{ \vr(T,\cdot) \varphi(T,\cdot)}- \intO{\vr_{0} \varphi(0,\cdot)}\\=
 \intTO{ \Big(\vr \partial_t \varphi + \vr \vu \cdot \nabla \varphi\Big)},
}
for all test functions $ \varphi \in C^1(\QTb)$;
\item{(ii)} $ \vr \vu \in C_w([0,T];L^{2}(\Omega,\R^3))$, and (\ref{sysSZb}) is satisfied in the weak sense
\eq{
&\intO{ (\vr \vu) (T,\cdot) \cdot \bm{\psi}(T,\cdot)}- \intO{ \vc{m}_0 \cdot \bm{\psi}(0,\cdot)}\\
&\quad=
\intTOB{\vr \vu \cdot \partial_ t \bm{\psi}  + \vr \vu \otimes \vu : \nabla \bm{ \psi} } \\
&\qquad+ \intTOB{\pi_\ep(Z) \Div \bm{\psi} +p(Z) \Div \bm{\psi} - \vS(\vu) : \Grad\bm{ \psi}   },
}
for all test functions $\bm{\psi} \in C_c^1([0,T] \times \Omega,\R^3)$;
\item{(iii)}
$ Z \in C_w([0,T];L^\infty(\Omega))$,  and \eqref{sysSZc} is satisfied in the weak sense
\eq{ 
\intO{Z(T,\cdot) \varphi(T,\cdot) \, \dx- \int_\Omega Z_{0} \varphi(0,\cdot)}\\=
 \intTO{\Big(Z \partial_t \varphi + Z \vu \cdot \nabla \varphi\Big)}, 
}
for all test functions $ \varphi \in C^1(\QTb)$;
\item{(iv)}
the energy inequality
\eq{
 {\cal E}(\vr,\vu, Z)(T)  + \intTO{ \Big(\mu |\nabla \vu|^2 + (\mu +\lambda)(\Div \vu)^2\Big)}\\ \leq {\cal E}(\vr_{0},\vu_{0},Z_0)
}
holds for a.a $T>0$, where
\begin{equation}
{\cal E}(\vr,\vu,Z) = \intO{\Big(\frac 12 \vr |\vu|^2 +Z\Gamma(Z)\Big)},
\end{equation}
\eq{
\Gamma(Z)=\int_0^Z\frac{\pi(s)+p(s)}{s^2}\, {\rm d}s.}
\end{description}
\end{df}

Our second goal is to show the convergence of the weak solutions to the system (\ref{sysSZ}), to the solutions of the limit system:
\begin{subequations}\label{sysSZ_lim}
\begin{equation}\label{cont1_lim}
 \pt\vr+ \Div (\vr\vu) = 0,
\end{equation}
\begin{equation}\label{mom1_lim}
\partial_t (\vr\vu) + \Div (\vr\vu \otimes \vu) + \Grad\pi +\Grad p(Z)-\Div\vc{S}(\vu) = \vc{0},
\end{equation}
\begin{equation}\label{Z1_lim}
\pt Z+\Div(Z\vu)=0,
\end{equation}
\begin{equation}\label{cons_lim}
0\leq Z\leq 1,\quad c_\star\vr\leq Z\leq c^\star\vr
\end{equation}
\begin{equation}\label{div_free}
\Div\vu=0\ \text{in} \ \{Z=1\}
\end{equation}
\begin{equation}\label{rel_pi}
\pi\geq0\ \text{in} \ \{Z=1\},\quad \pi=0\ \text{in} \ \{Z<1\},
\end{equation}
\end{subequations}
In the limit $\ep$, the uncongested (compressible) flow changes to the incompressible when the density $\vr$ hits the value $\vrs$. When this is the case the dynamics of the system is modified abruptly, meaning that the transition from the uncongested motion $(\vr<\vrs)$ to the congested motion $(\vr=\vrs)$ is very sudden.

The weak solutions to the limit system are defined as follows:
\begin{df}[Weak solution of the limit system with $Z$]\label{Def2}
A quadruple $(\vr,\vu,Z,\pi)$ is called a weak solution to \eqref{sysSZ_lim} with \eqref{boundary} and \eqref{initialZ} if equations \eqref{cont1_lim}, \eqref{mom1_lim}, \eqref{Z1_lim}
are satisfied in the sense of distributions, the constraints \eqref{cons_lim} satisfied a.e. in $(0,T)\times \Omega$, the divergence free condition \eqref{div_free} is satisfied a.e. in $\{Z=1\}$, 
and the following regularity properties hold
\begin{align*}
&\vr\in C_w([0,T];L^\infty(\Omega)), \\
&Z\in C_w([0,T];L^\infty(\Omega)),  \\
& \vu \in L^2(0,T;W^{1,2}_0(\Omega, \R^3)),\quad \vr|\vu|^2 \in L^{\infty}(0,T; L^1(\Omega)),\\
&\pi\in {\cal M}^+ ((0,T)\times \Omega).
\end{align*}
Moreover, $\pi$ is sufficiently regular so that the condition
\eqh{(1-Z)\pi=0,}
is satisfied in the sense of distributions on $(0,T)\times \Omega$.
\end{df}
The convergence result reads as follows.
\begin{thm}\label{Th:lim}
Let $(\vr_{\ep},\vu_{\ep}, Z_{\ep})_{\{\ep>0\}}$ be a a sequence of weak solutions to the approximate system  (\ref{sysSZ}), \eqref{pie}. Then, for $\ep\to 0$, the sequence $(\vr_{\ep},\vu_{\ep}, Z_{\ep})$ converges  to the weak solution of \eqref{sysSZ_lim} in the sense of Definition \ref{Def2}

More precisely,
\eqh{
\vr_{\ep}\rightarrow \vr \quad  \text{in }\ C_{w}([0,T];L^{\infty}(\Omega)), \quad \text{and weakly in } L^p((0,T)\times\Omega),\\
Z_{\ep}\rightarrow Z \quad \text{in }\ C_{w}([0,T];L^{\infty}(\Omega)) , \quad \text{and strongly in } L^p((0,T)\times\Omega),
}
for any $p<\infty$, and
\eqh{ \vu_{\ep} \rightarrow \vu \qquad \text{weakly in }  L^2(0,T;W^{1,2}(\Omega, \R^3)).}
Moreover,
\eqh{
&\pi_{\ep}(Z_\ep) \longrightarrow \pi  \quad\text{weakly\ in }\quad {\cal M}^+ ((0,T)\times \Omega).
}
\end{thm}

The paper is organized in the following manner. In Section \ref{Sect:exist},
we present details of approximation and prove the existence of solutions to the system \eqref{sysSZ} for $\ep$ fixed. Then, in Section \ref{Sec:lim}, we recover the two-phase system \eqref{sysSZ_lim} by letting $\ep\to 0$. After this, in Section \ref{Sec:rec} we recover the solution to the original two-phase system \eqref{sysSl}. Finally, in Section \ref{Sec:num} we briefly describe the numerical scheme and present computational examples that illustrate the behaviour of approximate solutions to the system \eqref{sysSl}.

\section{The existence of solution for $\ep$ fixed}\label{Sect:exist}
When $\ep$ is fixed, say $\ep=1$, the system \eqref{sysSZ} shares the features of the system considered in \cite{PeZa} and in the recent paper \cite{M3NPZ}; in our proof of existence of solutions we will  recall some elements of these two approaches in order to avoid repetitions. 

\subsection{Formulation of the approximate problem}
The first level of approximation introduces  truncation parameter $\delta$ in the singular pressure  and artificial pressure $\kappa\vr^K$, with $K$ sufficiently large to be determined in the course of the proof. We consider
\begin{subequations}\label{sysSZd}
\begin{equation}\label{cont1}
 \pt\vrd+ \Div (\vrd\vud) = 0,
\end{equation}
\begin{equation}\label{mom1}
\partial_t (\vrd\vud) + \Div (\vrd\vud \otimes \vud) + \Grad\pi_{\delta}(Z_\delta) +\Grad p_\kappa(Z_\delta)-\Div\vc{S}(\vud) = \vc{0},
\end{equation}
\begin{equation}\label{Z1}
\pt Z_\delta+\Div(Z_\delta\vud)=0,
\end{equation}
\end{subequations}
with $\kappa,\ \delta>0$, and $\pi_{\delta}$, $p_\kappa$ given by
\eq{\label{pi:ap}
\pi_{\delta}(Z_\delta)=\left\{
\begin{array}{lll}
\frac{Z_\delta^\alpha}{(1-Z_\delta)^\beta}&\mbox{if}& Z_\delta<1-\delta,\\
\frac{Z_\delta^\alpha}{\delta^\beta}&\mbox{if}& Z_\delta\geq1-\delta,
\end{array}
\right.
}
\eqh{
p_\kappa(Z_\delta)=\kappa Z_\delta^K+Z_\delta^\gamma.
}
We drop the subindex $\delta$ when no confusion can arise, and  we introduce the notion of weak solution for the system \eqref{sysSZd}.

\begin{df}[Weak solution of the approximate system ]\label{weaksolutionaux}
Suppose that the initial conditions satisfy (\ref{initialdata2}).
We say that the triplet  $(\vr,\vu,Z)$ is a weak solution to the problem \eqref{sysSZd} with the initial and boundary conditions (\ref{initialZ}) and (\ref{boundary}) if
\begin{equation*}
(\vr,\vu, Z) \in  L^\infty(0,T;L^K(\Omega))\times L^2(0,T;W^{1,2}_0(\Omega,\R^3)) \times L^\infty(0,T; L^K(\Omega)) ,
\end{equation*}
and for any $T>0$ we have:
\begin{description}
\item{(i)}
$ \vr \in C_w([0,T];L^K(\Omega))$, and \eqref{cont1} is satisfied in the weak sense
\eq{ 
\intO{ \vr(T,\cdot) \varphi(T,\cdot)}- \intO{\vr_{0} \varphi(0,\cdot)}\\=
 \intTO{ \Big(\vr \partial_t \varphi + \vr \vu \cdot \nabla \varphi\Big)},
}
for all test functions $ \varphi \in C^1(\QTb)$;
\item{(ii)} $ \vr \vu \in C_w([0,T];L^{\frac{2K}{K+1}}(\Omega,\R^3))$, and (\ref{mom1}) is satisfied in the weak sense
\eq{\label{weak_mom}
&\intO{ (\vr \vu) (T,\cdot) \cdot \bm{\psi}(T,\cdot)}- \intO{ \vc{m}_0 \cdot \bm{\psi}(0,\cdot)}\\
&\quad=
\intTOB{\vr \vu \cdot \partial_ t \bm{\psi}  + \vr \vu \otimes \vu : \nabla \bm{ \psi} } \\
&\qquad+ \intTOB{\pi_{\delta}(Z) \Div \bm{\psi}+p_{\kappa}(Z) \Div \bm{\psi} - \vS(\vu) : \Grad\bm{ \psi}   },
}
for all test functions $\bm{\psi} \in C_c^1([0,T] \times \Omega,\R^3)$;
\item{(iii)}
$ Z \in C_w([0,T];L^K(\Omega))$, and \eqref{Z1} is satisfied in the weak sense
\eq{
\intO{Z(T,\cdot) \varphi(T,\cdot) \, \dx- \int_\Omega Z_{0} \varphi(0,\cdot)}\\=
 \intTO{\Big(Z \partial_t \varphi + Z \vu \cdot \nabla \varphi\Big)}, 
}
for all test functions $ \varphi \in C^1(\QTb)$;
\item{(iv)}
the energy inequality
\eq{
 {\cal E}(\vr,\vu, Z)(T)  + \intTO{ \Big(\mu |\nabla \vu|^2 + (\mu +\lambda)(\Div \vu)^2\Big)}\\ \leq {\cal E}(\vr_{0},\vu_{0},Z_0)
}
holds for a.a. $T>0$, where
\begin{equation}
{\cal E}(\vr,\vu,Z) = \intO{\Big(\frac 12 \vr |\vu|^2 +Z\Gamma_{\kappa,\delta}(Z)\Big)},
\end{equation}
\eq{\label{defG}
\Gamma_{\kappa,\delta}(Z)=\int_0^Z\frac{\pi_{\delta}(s)+p_\kappa(s)}{s^2}\, {\rm d}s.}
\end{description}

\end{df}

We have the following existence result for solutions defined by Definition \ref{weaksolutionaux} (see also \cite{M3NPZ}, Theorem 2).
\begin{thm}\label{mainthm2}
Let $\vS$ satisfy (\ref{S}), $K>6,\ \beta>5/2$, $\alpha\geq0$, $\kappa,\delta,\ep$ be fixed and positive, and the initial data $(\vr_0,\vc{m}_0,Z_0)$ satisfy (\ref{initialdata2}).

Then there exists a weak solution $(\vr,\vu, Z)$ to problem \eqref{sysSZd}, \eqref{pi:ap} with boundary conditions \eqref{boundary}, in the sense of Definition \ref{weaksolutionaux}. 

Moreover, $(Z,\vu)$ solves \eqref{Z1}  in the renormalized sense, i.e. 
$(Z,\vu)$, extended by zero outside of $\Omega$, satisfies
\begin{equation}\label{renorent}
\partial_t b(Z) + \Div(b(Z)\vu) + \big(b'(Z)Z-b(Z)\big) \Div \vu = 0,
\end{equation}
in the sense of distributions on $(0,T)\times \R^3$, 
where
\begin{equation}\label{regb}
b \in C^1(\R), \quad b'(z)=0,\quad \forall z \in \R~\text{large enough}.
\end{equation}
In addition,
\eq{\label{cZc}
0 \leq c_\star \vr \leq Z \leq c^\star \vr\quad
\text{a.e. in }(0,T)\times \Omega.}
\end{thm}

The proof of  Theorem \ref{mainthm2}  requires further modification of the system \eqref{sysSZd} with two additional approximation levels involving the parabolic regularisation of  two continuity equations and the Galerkin approximation of the velocity. This approximation allows, in particular, to deduce inequalities \eqref{cZc}. At this point existence of regular solution can be obtained following the arguments presented in \cite{M3NPZ}. Also the compactness arguments needed to recover system \eqref{sysSZd} are analogous, therefore we skip this part and focus  only on the a-priori estimates needed to perform the limit passages $\delta\to0$, $\kappa\to 0$.

Having the existence of solutions to the system (\ref{sysSZd})--\eqref{pi:ap}, we show that this solution can be used to recover the weak solution to the system \eqref{sysSZ}, where only the parameter $\ep$ is present. 

\begin{thm}\label{Th:ep}
Let $\ep$ be fixed and let $(\vr_{\kappa,\delta},\vu_{\kappa,\delta}, Z_{\kappa,\delta})$ be a weak solution to the approximate system  (\ref{sysSZd})--\eqref{pi:ap} established in Theorem \ref{mainthm2}. Then, for $\delta,\kappa\to 0$, the sequence $(\vr_{\kappa,\delta},\vu_{\kappa,\delta}, Z_{\kappa,\delta})_{\kappa,\delta>0}$ converges to the weak solution of \eqref{sysSZ} in the sense of Definition \ref{weaksolutionaux}.

More precisely,
\eqh{
\vr_{\kappa,\delta}\rightarrow \vr \quad  \text{in }\ C_{w}([0,T];L^{\infty}(\Omega)), \quad \text{and weakly in } L^p((0,T)\times\Omega),\\
Z_{\kappa,\delta}\rightarrow Z \quad \text{in }\ C_{w}([0,T];L^{\infty}(\Omega)) , \quad \text{and strongly in } L^p((0,T)\times\Omega),
}
for any $p<\infty$, and
\eqh{ \vu_{\kappa,\delta} \rightarrow \vu \qquad \text{weakly in }  L^2(0,T;W^{1,2}(\Omega, \R^3)).}
Moreover,
\eqh{
&\pi_{\delta}(Z_{\kappa,\delta}) \longrightarrow \pi_\ep(Z)  \quad\text{strongly\ in }\quad L^1((0,T)\times \Omega).
}
\end{thm}

In the next two sections we prove Theorem \ref{Th:ep}. Starting from the weak solution to the system \eqref{sysSZd}, we first derive  uniform bounds and then we let $\delta\to 0$. The modification  of the reasoning needed to perform the limit passage $\kappa\to 0$ is explained at the end.

\subsection{Uniform estimates}
In this section we obtain the uniform estimates for the weak solutions to the system \eqref{sysSZd}--\eqref{pi:ap}. This involves three a-priori estimates: the standard energy estimate, and two estimates involving the application of the Bogovskii operator, one of which gives us the uniform bound for the singular pressure and the other the uniform integrability of the pressure. All of the aforementioned estimates are essential for performing the limit passage $\delta,\kappa\to 0$, as well as the last limit passage $\ep\to0$. However, as we shall see later on, the uniform integrability of the pressure will no longer be valid in this case.

\subsubsection{The energy estimate}
We present a formal computation that can be made rigorous at the level of the Galerkin approximation of the velocity. Multiplying the momentum equation \eqref{mom1} by $\vu$ and integrating by parts with respect to space, yield the energy equality
\eq{\label{en:test}
\Dt \intO{\frac 12 \vr |\vu|^2}+\intO{\Grad\lr{\pi_{\delta}(Z)+p_\kappa(Z)}\cdot\vu}+\intO{\vS(\vu):\Grad\vu}=0.}
Note that the second term in our case is different than in \cite{PeZa}, there is no additional $\vr^*$ in front of the gradient. This however allows us to proceed straightforwardly, for reeder's convenience we repeat here the derivation
\eqh{
&\intO{\Grad\lr{\pi_{\delta}(Z)+p_\kappa(Z)}\cdot\vu}=\intO{\frac{\pi'_{\delta}(Z)+p'_\kappa(Z)}{Z}\Grad Z\cdot(Z\vu)}\\
&=-\intO{Q_{\kappa,\delta}(Z)\Div(Z\vu)}=\intO{Q_{\kappa,\delta}(Z)\pt Z}= \Dt\intO{Z\Gamma_{\kappa,\delta}(Z)}
}
where we used the equation \eqref{Z1}, we denoted $Q'_{\kappa,\delta}(Z)=\frac{\pi'_{\delta}(Z)+p'_\kappa(Z)}{Z}$, and $\Gamma_{\kappa,\delta}$ is a solution of the following ODE:
\eqh{\Gamma_{\kappa,\delta}(Z)+Z\Gamma'_{\kappa,\delta}(Z)=Q_{\kappa,\delta}(Z).}
Using the definition of $Q_{\kappa,\delta}$ and $\pi_{\kappa,\delta}$ we can express $\Gamma_{\kappa,\delta}$ as in \eqref{defG}.
Integrating \eqref{en:test} with respect to time, and using the definition of the stress tensor \eqref{S} we get the following uniform estimates
\begin{equation} \label{un_delta}
\begin{gathered}
\sup_{t\in[0,T]}\lr{\|\sqrt{\vr_\delta} \vu_\delta (t)\|_{L^2(\Omega)} 
+\left\|Z_\delta \Gamma_{\kappa,\delta}(Z_\delta)(t) \right\|_{L^1(\Omega)} }\leq C,\\
\intT{{\|\vu_\delta\|_{W^{1,2}(\Omega,\R^3)}^2}} \leq C.
\end{gathered}
\end{equation}
\subsubsection{The integrability of the pressure}
The energy estimate obtained above is insufficient to control the  $L^1$ norm of the pressure, because the singularity appearing in $\Gamma_{\kappa,\delta}(Z)$ at $Z=1$ is of lower order than the singularity for $\pi_{\kappa,\delta}(Z)$ (take f.i. $\alpha=2$, $\beta=1$ in \eqref{pie} and use \eqref{defG}). Therefore, further estimates are needed. The first of them is obtained by testing the momentum equation by the function
\eq{\label{testB}
\bm{\psi}=\phi(t){\cal B}\lr{Z-\frac{1}{|\Omega|}\int_\Omega{Z(t,y)\,{\rm d}y}},}
where $\phi$ is smooth and compactly supported in the interval  $(0,T)$, and ${\cal B}$ is the  Bogovskii operator, whose main properties can be found, for example, in \cite[Lemma 3.17]{NS}, and  \cite[Appendix]{PeZa}. In particular, ${\cal B}$ is a linear operator such that for any $f\in L^p(\Omega)$, $1<p<\infty$, $\int_\Omega{f}(y){\rm d}y=0$
\[\Div\mathcal{B}(f)=f \quad \text{a.e. in } \Omega,\]
and
\[\|\nabla \mathcal{B}(f)\|_{L^p(\Omega, \R^3\times\R^3)} \leq c(p,\Omega)\|f\|_{L^p(\Omega)}, \quad 1<p<\infty.\]
Moreover, if $f=\Div g$, with  $g\in L^q(\Omega,\R^3)$, $\Div g\in L^p(\Omega)$, $1<q<\infty$, then
\[\|\mathcal{B}(f)\|_{L^q(\Omega,\R^3)}\leq c(q,\Omega)\|g\|_{L^q(\Omega,\R^3)}.\]
\begin{rmk}
Note that, alike for $\vr$, the integral of $Z$ over the space is an invariant  of the motion, therefore we have
\eq{\intO{Z(t,x)}=\intO{Z_0(x)}<|\Omega|,}
according to \eqref{initialdata2}, at least for sufficiently regular solutions. In what follows we denote $\intO{Z_0(x)}=M_Z$.
\end{rmk}

\begin{rmk}
Function $\varphi$ must be of a certain regularity in order to use it as a test function in \eqref{weak_mom}. In fact, it follows from the construction of the solution in \cite{M3NPZ}, that taking $K$ sufficiently large guarantees admissibility of this function.
\end{rmk}

Using in \eqref{weak_mom} the test function \eqref{testB} results in the following equality
\eqh{
&\intTO{\phi\lr{\pi_{\delta}(Z_\delta)+p_{\kappa}(Z_\delta)} \lr{Z_\delta-\frac{1}{|\Omega|}\int_\Omega{Z_\delta\,{\rm d}y}}}\\
&\quad=
-\intTO{\vrd \vud \cdot \partial_ t \bm{\psi}}  -\intTO{ \vrd \vud \otimes \vud : \nabla \bm{ \psi} } \\
&\qquad+ \intTO{\vS(\vud) : \Grad\bm{ \psi}  },
}
whose r.h.s. can be bounded using the uniform estimates \eqref{un_delta} together with \eqref{cZc}, provided that $K$ is sufficiently large, say $K>4$.
Therefore, one gets the following estimate, which is now uniform with respect to $\kappa$ and $\delta$
\eq{\label{Bog}
\intTO{\phi\lr{\pi_{\delta}(Z_\delta)+p_{\kappa}(Z_\delta)}  \lr{Z_\delta-\frac{1}{|\Omega|}\int_\Omega{Z_\delta\,{\rm d}y}}}\leq C.
}
We then consider two complementary subsets of $(0,T)\times \Omega$: $\Sigma_1=\{Z_\delta(t,x)<Z^*\}$ and $\Sigma_2=\{ Z_\delta(t,x)\geq Z^*\}$, for $\frac{M_Z}{|\Omega|}<Z^*<1$. The l.h.s. of \eqref{Bog} can be easily controlled on $\Sigma_1$, because on this subset $Z_\delta$ stays far away from the singularity, for $\Sigma_2$ we have
\eqh{\intTO{\phi\pi_{\delta}(Z_\delta) \lr{Z_\delta-\frac{1}{|\Omega|}\int_\Omega{Z_\delta\,{\rm d}y}} \vc{1}_{\Sigma_2}}\\
\geq 
\lr{Z^*-\frac{M_Z}{|\Omega|}}\intTO{\phi\pi_{\delta}(Z_\delta) \vc{1}_{\Sigma_2}},
}
and so, from \eqref{Bog} it follows that
\eq{\label{un_pi}
\|\pi_{\delta}(Z_\delta)\|_{L^1((0,T)\times\Omega)}+\|Z_\delta\pi_{\delta}(Z_\delta)\|_{L^1((0,T)\times\Omega)}\leq C,}
as well as
\eq{\kappa\|Z_\delta\|_{L^{K+1}((0,T)\times\Omega)}^{K+1}\leq C.
}
The estimate \eqref{Bog} implies that the $L^{K+1}$ norm of  $Z_\delta$ is bounded uniformly in $\delta$, but not uniformly in $\kappa$, which suggests that the passage to the limit $\delta\to0$ should be performed as first.

\subsubsection{The equi-integrability of the singular pressure}
Because $\pi_{\delta}$ is a nonlinear function of $Z$, identification of the limit in this term, after letting $\delta\to0$ cannot be justified only by the uniform $L^1$ bound. Following the idea from \cite{FePeRoSch}, we can prove 
that the pressure $\pi_{\delta}$ enjoys some additional estimate near to the singularity $Z=1$. Indeed,
for $K>6$, $\beta>5/2$ and 
\begin{equation}\label{eta_function}
\eta_\delta(s) = \begin{cases}  \,\,-\log (1-s) & \text{ if } s \leq 1-\delta, \\
                                 \,\,-\log (\delta)  & \text{ if } s > 1-\delta, \end{cases}
\end{equation} 
uniformly with respect to $\delta$ one has
\begin{equation}\label{equiintegr}
\intTO{\pi_{\delta}(Z)\eta_\delta(Z)} \leq C. 
\end{equation} 
The proof follows by testing the momentum equation \eqref{weak_mom} by the function of the form
\eq{\label{testBII}
\bm{\psi}=\phi(t){\cal B}\lr{\eta_\delta(Z_\delta)-\frac{1}{|\Omega|}\int_\Omega{\eta_\delta(Z_\delta)\,{\rm d}y}},}
where $\phi$ is smooth and compactly supported in the interval  $(0,T)$. For the details of this estimate we refer to \cite{PeZa}, for $\beta>3$ and to \cite{FeLuMa} for $\beta>\frac{5}{2}$. 
One of the difficulties in the proof of analogue of  \eqref{equiintegr} presented in \cite{PeZa} concerned the renormalization of the equation for $\frac{\vr}{\vrs}$. Here this problem does not appear anymore, since $(Z_\delta,\vud)$ is by definition a distributional solution of the renormalized transport equation \eqref{renorent}.
\subsection{Passage to the limit $\delta\to 0$}
\subsubsection{Convergences following from the uniform estimates}
Using the uniform estimates  \eqref{cZc}, \eqref{un_delta}, and the H\"older inequality, we can deduce that up to a subsequence
\eq{\label{conv}
&Z_\delta\rightarrow Z \qquad \text{weakly-*  in }  L^{\infty}(0,T;L^{K}(\Omega)), \\
&\vrd\rightarrow \vr \qquad \text{weakly-*  in }  L^{\infty}(0,T;L^{K}(\Omega)), \\
& \vud \rightarrow \vu \qquad \text{weakly in }  L^2(0,T;W^{1,2}(\Omega, \R^3)).
}
Using the continuity equation \eqref{cont1} and \eqref{Z1}, the the two first convergences can be strengthened to
\eq{
&Z_\delta\rightarrow Z \qquad \text{in }\ C_w([0,T];L^{K}(\Omega))\\
&\vrd\rightarrow \vr \qquad \text{in }\ C_w([0,T];L^{K}(\Omega)),
}
which together with the weak convergence of the velocity gradient, after using the momentum equation \eqref{mom1} implies that 
\eq{\vrd\vud \rightarrow \vr \vu \quad \text{in }\  C_w([0,T]; L^{\frac{2K}{K+1}}(\Omega,\R^3)).}
This in turn, assuming that $K$ is sufficiently large so that the imbedding of $L^{\frac{2K}{K+1}}(\Omega,\R^3)$ to $W^{-1,2}(\Omega,\R^3)$ is compact, implies that
\begin{equation}\label{cvg_conv}
\vrd\vud\otimes \vud \rightarrow \vr \vu \otimes \vu \quad \text{weakly in } L^q((0,T)\times \Omega, \R^3\times\R^3),
\end{equation}
for some $q>1$. Moreover, the uniform bound \eqref{un_pi} together with the growth condition $\beta>5/2$ in \eqref{pi:ap} imply that
\eq{\label{unZ}
Z\leq 1\quad \text{a.a.}\ (t,x)\in(0,T)\times\Omega,
}
and thus also
\eqh{
\vr\leq \frac{1}{c_\star}\quad \text{a.a.}  (t,x)\in(0,T)\times\Omega,
}
on the account of \eqref{cZc}.
Finally, from the uniform bounds \eqref{un_pi} we can extract the subsequences such that
\eq{
&\pi_{\delta}(Z_\delta) \rightarrow\Ov{\pi(Z)}  \quad\text{weakly\ in }\quad {\cal M}^+((0,T)\times\Omega),\\
&Z_\delta \pi_{\delta}(Z_\delta) \rightarrow\Ov{Z\pi(Z)}  \quad\text{weakly\ in }\quad {\cal M}^+((0,T)\times\Omega),
}
for some $\Ov{\pi(Z)}$, $\Ov{Z\pi(Z)} $ that need to be determined. Note, however, that \eqref{equiintegr} together with De~La~Vall\'ee-Poussin criterion  allow us to deduce the equi-integrability of the pressure, therefore the first limit can be strengthened to
\eq{\label{piL1}
&\pi_{\delta}(Z_\delta) \rightarrow\Ov{\pi(Z)}  \quad\text{weakly\ in }\quad L^1((0,T)\times \Omega).
}
As for the second convergence, we cannot say that much immediately. However, using the fact that $\pi_{\delta^1}(\cdot) \geq \pi_{\delta^2}(\cdot)$ provided $\delta^1\leq \delta^2$, we can estimate for any smooth, nonnegative,  compactly supported function $\phi(x,t)$
\eq{ \label{pi_pi}
&\liminf_{\delta\rightarrow 0} \intTO{\phi\lr{\pi_{\delta}(Z_\delta)Z_\delta-\overline{\pi(Z)}Z}}\\
&\geq \liminf_{\delta\rightarrow 0}\intTO{\phi\lr{\pi_{\delta^*}(Z_\delta)Z_\delta-\overline{\pi(Z)}Z}}\\
&\geq \liminf_{\delta\rightarrow 0}\intTO{\phi\lr{\Ov{\pi_{\delta^*}(Z)}-\overline{\pi(Z)}}Z}
}
where the last inequality is a consequence of the fact that $\cdot\mapsto\pi_{\delta^*}(\cdot)$ is non-decreasing function for fixed $\delta^*$. Using again equi-integrability of the pressure and letting $\delta^*\to 0$ we verify that the r.h.s. of \eqref{pi_pi} vanishes and thus
\eq{\label{one_side}
\Ov{\pi(Z)}\geq Z\Ov{\pi(Z)}}
in the sense of distributions.
In order to say something more, we need to investigate the strong convergence of the sequence $Z_\delta$, which is the purpose of the next section.
\subsubsection{Strong convergence of $Z_\delta$}
It was  observed in \cite{M3NPZ}, that the strong convergence of $Z_\delta$ does not imply the strong convergence of $\vr_\delta$.  The proof of the strong convergence of $Z_\delta$ requires an analogue of the effective flux equality for the barotropic compressible Navier-Stokes equations. In our case it can be written as follows.
\begin{lemma}\label{effectivestep1}
Let $\vr_\delta,\vu_\delta, Z_\delta$ be the sequence of approximate solutions enjoying the properties from above. Then, at least for a subsequence
\eq{ \label{evf}
\lim_{\delta \to 0^+}  \intTO{\psi  \phi \big( \pi_{\delta}(Z_\delta) +p_{\kappa}(Z_\delta)-(\lambda+2\mu) \Div \vu_\delta\big) Z_\delta}\\
=   \intTO{\psi \phi  \big( \overline{\pi(Z)}+\Ov{p_{\kappa}(Z)} -(\lambda+2\mu) \Div \vu\big) Z }
}
for any $ \psi \in C_c^\infty((0,T)) $ and $\phi \in C_c^\infty(\Omega)$.
\end{lemma}

The proof of this fact can be seen as a special case of an analogous result proven in \cite[cf. Lemma 11]{M3NPZ}. 
On the account of \eqref{one_side}  and the monotonicity of $p_\kappa(\cdot)$ we obtain from \eqref{evf} that
\eq{\label{in1}
\lim_{\delta \to 0^+}  \intTO{ \psi \phi Z_\delta\, \Div \vu_\delta }-\intTO{\psi \phi  Z \, \Div \vu  }\geq0,
}
for any $\psi\phi\geq0$.
Recall that $Z_\delta$ satisfies the renormalized continuity equation \eqref{renorent} with $b$ specified in \eqref{regb}. 
By density argument and standard approximation technique, we may extend the validity of (\ref{renorent}) to functions $b \in C([0,\infty) \cap C^1((0,\infty))$ such that
\[
|b'(t)| \leq C t^{-\lambda_0}, \qquad \lambda_0 <-1, \quad t \in (0,1],
\]
\[
|b'(t)| \leq C t^{\lambda_1}, \qquad -1<\lambda_1 \leq \frac{q}{2}-1, \quad t \geq 1.
\]
The renormalization technique applied to the barotropic Navier-Stokes system is due to DiPerna and Lions \cite{DPL}, and the above extension can be found, for example, in \cite{FN}. 

We can now write the renormalized continuity equation with $b(Z_\delta)=Z_\delta\log Z_\delta$:
\eqh{
\partial_t(Z_\delta\log Z_\delta) + \Div(Z_\delta\log Z_\delta\vud) = -Z_\delta\Div \vud \quad \text{in}\quad \mathcal{D}'((0,T)\times \mathbb{R}^3).
}
Passing to the limit $\delta\to 0^+$, we hence obtain
\eqh{\partial_t (\overline{Z \log Z} )+ \Div(\overline{Z\log Z}\vu) =-\overline{Z \Div \vu} \quad \text{in}\quad \mathcal{D}'((0,T)\times \mathbb{R}^3).}
Writing an analogous equation for the limit function $Z$ and subtracting it from the above, we obtain
\begin{equation}\label{ineg_log}
\partial_t(\overline{Z\log Z}-Z\log Z) + \Div[(\overline{Z\log Z}-Z\log Z)\vu]=Z\Div \vu-\overline{Z\Div \vu},
\end{equation}
satisfied in the sense of distributions on $(0,T)\times \mathbb{R}^3$. Integrating the above equation with respect to time and space, and using the convexity of the function $s\mapsto s\log s$ we get from \eqref{ineg_log} that
\[\intTO{\overline{Z\Div \vu}} \leq \intTO{Z\Div \vu},\]
which is an opposite to \eqref{in1}. Therefore, recalling \eqref{ineg_log} we see that $\overline{Z\log Z}(t,x)=Z\log Z(t,x)$ almost everywhere in $(t,x)\in(0,T)\times\Omega$, which yields the strong convergence of $Z_\delta$ in $L^p((0,T)\times\Omega)$ for any $p<K+1$. With this at hand, we verify that \eqref{piL1} can be replaced by
\eqh{
&\pi_{\delta}(Z_\delta) \rightarrow \pi(Z)  \quad\text{strongly\ in }\  L^1((0,T)\times \Omega),
}
similarly
\eqh{
&p_{\kappa}(Z_\delta) \rightarrow p_\kappa(Z)  \quad\text{strongly\ in }\  L^1((0,T)\times \Omega).
}
In order to complete the proof of Theorem \ref{Th:ep}, we have to show that we are allowed to let $\kappa\to 0$. Note that all the uniform estimates obtained above stay in force independently of $\kappa$. Indeed, the only issue here is to see that the uniform $L^{K+1}$ bound on $Z_\kappa$, that was previously obtained from \eqref{un_pi} follows now directly from \eqref{unZ}. On the account of \eqref{cZc}, the same is true also for the sequence $\vr_\kappa$. The proof of Theorem \ref{Th:ep} is now complete. $\Box$

\section{Passage to the limit $\ep\to 0$} \label{Sec:lim}
The purpose of this section is to prove our main Theorem \ref{Th:main}. For technical reasons we first perform the limit $\ep\to0$ in the auxiliary system \eqref{sysSZ} proving Theorem \ref{Th:lim} and then we prove equivalence between systems \eqref{sysSZ_lim} and \eqref{sysSl} in the certain class of solutions.

\subsection{Convergence following from the uniform estimates}
The estimates performed in the previous section give rise to several estimates that are uniform with respect to $\ep$. Indeed, passing to the limit in the energy estimate we obtain
\begin{equation} \label{un_ep}
\begin{gathered}
\sup_{t\in[0,T]}\lr{\|\sqrt{\vr_\ep} \vu_\ep (t)\|_{L^2(\Omega)} 
+\left\|Z_\ep \Gamma_{\ep}(Z_\ep)(t) \right\|_{L^1(\Omega)} }\leq C,\\
\intT{{\|\vu_\ep\|_{W^{1,2}(\Omega,\R^3)}^2}} \leq C.
\end{gathered}
\end{equation}
Passing to the limit in \eqref{unZ}  and in \eqref{cZc} we obtain 
\eq{\label{Zr_ep}
0\leq Z_\ep\leq1,\quad 0\leq c_\star \vr_\ep\leq Z_\ep\leq c^\star\vr_\ep,}
in particular both sequences $Z_\ep$, $\vr_\ep$ are uniformly bounded in $L^p((0,T)\times\Omega)$ for $p\leq \infty$. Therefore, by means of  the arguments from the previous section we get, up to the subsequence, that
\eq{\label{conv_ep}
 &\vu_\ep \rightarrow \vu \qquad \text{weakly in }  L^2(0,T;W^{1,2}(\Omega, \R^3),\\
&Z_\ep\rightarrow Z \qquad \text{in }\ C_{w}([0,T];L^{\infty}(\Omega)),\\
&\vr_\ep\rightarrow \vr \qquad \text{in }\ C_{w}([0,T];L^{\infty}(\Omega)).
}
This information allows us to pass to the limit in all terms of the system \eqref{sysSZ}, apart from the nonlinear pressure terms $p(Z)$ and $\pi_\ep(Z)$. Repeating the Bogovskii type estimate with the test function \eqref{testB} we obtain
\eq{\label{un_pi_ep}
\|Z_\ep p(Z_\ep)\|_{L^1((0,T)\times\Omega)}+\|\pi_{\ep}(Z_\ep)\|_{L^1((0,T)\times\Omega)}+\|Z_\ep\pi_{\ep}(Z_\ep)\|_{L^1((0,T)\times\Omega)}\leq C,}
however, the estimate \eqref{equiintegr} does not hold anymore. Therefore the convergence in the sense of measures is the most we can hope for, we have
\eq{\label{lim_pi_ep}
&\pi_{\ep}(Z_\ep) \rightarrow\pi  \quad\text{weakly\ in }\quad {\cal M}^+((0,T)\times\Omega),\\
&Z_\ep \pi_{\ep}(Z_\ep) \rightarrow\pi_1  \quad\text{weakly\ in }\quad {\cal M}^+((0,T)\times\Omega).
}
For the background pressure, due to \eqref{Zr_ep}, we have
\eqh{
p(Z_\ep)\to \Ov{p(Z)}\quad\text{weakly\ in }\quad L^{p}((0,T)\times\Omega),
}
for any $p<\infty$.
At this point, we can identify the second limit in \eqref{lim_pi_ep} using the explicit form of the pressure \eqref{pie}. We have
\eq{
Z_\ep\pi_\ep(Z_\ep)=\ep\frac{1}{(1-Z_\ep)^\beta}=\pi_\ep(Z_\ep)-\ep\frac{1}{(1-Z_\ep)^{\beta-1}},
}
thus letting $\ep\to0$ and observing that the last term converges to zero strongly, we obtain the relation
\eq{\label{step1}
\pi_1 =\pi  }
in the sense of the measures. The recovery of the constraint condition $(1-Z)\pi=0$ and the identification of the limit $\Ov{p(Z)}=p(Z)$ require stronger information about convergence of $Z_\ep$.

\subsection{Strong convergence of $Z_\ep$}
The nowadays well known technique of proving the strong convergence of the density in the compressible barotropic Navier-Stokes equations involves the study of propagation of the so called oscillation defect measure \cite{EF2001}. Such level of precision is not needed in our case, because, due to the singularity of the pressure argument $Z$, we have sufficiently high integrability of the pressure in order to apply the DiPerna-Lions technique \cite{DPL}. Nevertheless, a variant of effective viscous flux equality is still needed. It can be derived the same way as in Lemma \ref{effectivestep1}, after observing that the inverse divergence operator $\Grad\lap^{-1}[1_\Omega Z]$ is regular enough to be used as a test function in the limiting momentum equation. 
 
 With this information the statement of Lemma \ref{effectivestep1} adapted to the $\ep$-labelled sequences gives rise to the equality
\eq{ \label{evf_ep}
\lim_{\ep \to 0^+}  \intTO{\psi  \phi \big( \pi_{\ep}(Z_\ep)+p(Z_\ep) -(\lambda+2\mu) \Div \vu_\ep\big) Z_\ep}\\
=   \intTO{\psi \phi  \big( \pi +\Ov{p(Z)}-(\lambda+2\mu) \Div \vu\big) Z },
}
for any $ \psi \in C_c^\infty((0,T)) $ and $\phi \in C_c^\infty(\Omega)$.
 
Let us now explain the meaning of the product $Z \pi $ on the r.h.s. of \eqref{evf_ep}. To this end, one needs to come back to the limiting momentum equation
 $$ \partial_t (\vr\vu) + \Div (\vr\vu \otimes \vu) + \Grad\pi +\Grad \Ov{p(Z)}-\Div\vc{S}(\vu) = \vc{0},$$
 and use the bounds \eqref{un_ep}, \eqref{Zr_ep} to justify that $\pi$ is in fact more regular than it follows just from \eqref{lim_pi_ep}. Indeed, we have
 \begin{equation}\label{regu_pi}
\pi \in W^{-1,\infty}(0,T;W^{1,2}(\Omega)) \cup L^p(0,T;L^q(\Omega))\quad p,q>1.\end{equation}
Moreover from the equation for $Z$, we easily get
\eq{\label{regu_z}
Z \in C_w([0,T];L^\infty(\Omega))\cap C^1([0,T];W^{-1,2}(\Omega)).}
Regularizing in space and time the limits $Z$ and $\pi$ using the standard multipliers $\omega_n$,
$Z_n=Z\ast \omega_n$, $\pi_n=\pi\ast \omega_n$, we can clarify the meaning of $Z \pi$ by writing
\eq{\label{product}
Z\pi= Z_n\pi_n+(Z-Z_n)\pi_n+Z(\pi-\pi_n),}
and by passing to the limit with the support of mollifying kernel, see \cite{PeZa} for more details.

From \eqref{evf_ep} it follows that
\eq{
&(\lambda+2\mu)\intTO{\psi\phi\lr{\Ov{Z\Div\vu}-Z\Div\vu}}\\
&=\intTO{\psi\phi\lr{\pi_1-Z\pi}}+\intTO{\psi\phi\lr{\Ov{Zp(Z)}-Z\Ov{p(Z)}}}\\
&\geq\intTO{\psi\phi\lr{1-Z}\pi}\geq 0,
}
where to get to the r.h.s. of the above we have used subsequently: monotonicity of $p(\cdot)$, \eqref{step1}, and the limit of \eqref{Zr_ep}. Since both pairs $(Z_\ep,\vu_\ep)$ and $(Z,\vu)$ satisfy the renormalized continuity equation, we can use the renormalization in the form $b(z)=z\log z$ to justify that 
$$Z_\ep\to Z\quad \text{strongly in } L^p((0,T)\times \Omega),\quad \forall p<\infty.$$
Note however,  that this property is not transferred to the sequence $\vr_\ep$, for which we only have \eqref{conv_ep}. Nevertheless, using this information and formula \eqref{product} we can justify that
$$\pi_1=Z\pi,$$
which together with \eqref{step1} implies \eqref{rel_pi}. 

It remains to show the condition \eqref{div_free}, or rather its compatibility with the other conditions in system \eqref{sysSZ_lim}. This follows from the following lemma proven by Lions and Masmoudi in \cite{LM99}, that we recall here without the proof.
\begin{lemma}[Lemma 2.1, \cite{LM99}]
Let $\vu\in L^2(0,T; W^{1,2}_0(\Omega,\R^3))$ and $f\in L^2((0,T)\times\Omega)$ such that
$$\pt f+\Div(f \vu)=0\quad in \ (0,T)\times\Omega,\quad f(0)=f_0,$$
then the following two assertions are equivalent
\begin{itemize}
\item[(i)] $\Div\vu=0$~ a.e. on ~$\{f=1\}$ ~and  ~$0 \leq f_0\leq 1$,
\item[(ii)] $0\leq f(t,x)\leq1$.
\end{itemize}
\end{lemma}
Applying this lemma for $f=Z$ we conclude the proof of Theorem \ref{Th:lim}. $\Box$

\section{Recovery of the original system}\label{Sec:rec}
In Section \ref{Sec:lim} we proved that system \eqref{sysSZ_lim} possesses a weak solution in the sense of Definition \ref{Def2}. Our next aim is to prove that this solution can be identified with the solution to the original problem \eqref{sysSl}.
In other words, we need to deduce the existence of $\vr^*\in L^\infty((0,T)\times\Omega)$ satisfying the transport equation, such that the measure in the momentum equation vanishes for $\vr=\vr^*$.

We proceed similarly to  \cite{M3NPZ}. Note that 
\eqh{\frac{\vr_0}{Z_0}\Big|_{\{\vr_0=0\}} =\frac{\vr_0}{Z_0}\Big|_{\{Z_0=0\}}=\tilde{\vrs}>0.}
When extended by $0$ outside $\Omega$ the couples  $(\vr,\vu)$ and $(Z,\vu)$ satisfy the renormalized continuity equations, due to uniform $L^\infty$ bounds for both $\vr_\ep$ and $Z_\ep$. Therefore, we may test the equations \eqref{cont1_lim}, \eqref{Z1_lim} by $\omega_n(x-\cdot)$, where $\omega_n$ is a standard mollifier, which leads to
\begin{equation}\label{vr_mol}
\partial_t\vr_n+\Div(\vr_n\vu)= r^1_n,
\end{equation}
\begin{equation} \label{Z_mol}
\partial_t Z_n+\Div(Z_n\vu)= r^2_n,
\end{equation}
satisfied a.e. in $(0,T)\times\R^3$, where by $a_n$ we denoted $a\ast\omega_n$. It follows from the Friedrichs commutator lemma, see e.g. \cite[Lemma 10.12]{FN}, that $r^1_n$, and $r^2_n$ converge to $0$ strongly in $L^1((0,T)\times\R^3)$ as $n\to \infty$.

We now multiply \eqref{vr_mol}  by $\frac{1}{Z_n+\lambda}$, and \eqref{Z_mol} by $-\frac{\vr_n+\lambda\tilde{\vrs}}{(Z_n+\lambda)^2}$, with $\lambda>0$, and obtain, after some algebraic transformations, that
\eqh{\partial_t \left( \frac{ \vr_n + \lambda \tilde{\vrs}}{Z_n + \lambda } \right) + \Div \left[ \left( \frac{ \vr_n + \lambda \tilde{\vrs}}{Z_n + \lambda } \right) \vu \right]
-\left[ \frac{( \vr_n + \lambda \tilde{\vrs}) Z_n }{(Z_n + \lambda )^2} + \frac{\lambda \tilde{\vrs}}{ Z_n + \lambda  }\right] \Div \vu\\
= r^1_n \frac{1}{Z_n + \lambda}-r^2_n \frac{ \vr_n + \lambda \tilde{\vrs}}{(Z_n + \lambda)^2}.
}
By passing with $n\to \infty$, we get
\eq{\label{fin}
\partial_t \left( \frac{ \vr + \lambda \tilde{\vrs}}{Z+ \lambda } \right) + \Div \left[ \left( \frac{ \vr+ \lambda \tilde{\vrs}}{Z + \lambda } \right) \vu \right]
-\left[ \frac{( \vr + \lambda \tilde{\vrs}) Z }{(Z+ \lambda )^2} + \frac{\lambda \tilde{\vrs}}{ Z + \lambda  }\right] \Div \vu
= 0 .
}
We distinguish two cases:\\

{\it Case 1. } For $Z=0$, from \eqref{cons_lim} it follows that $\vr=0$ and therefore  $\frac{ \vr + \lambda \tilde{\vrs}}{Z+ \lambda }=\tilde{\vrs}$, and  $\frac{( \vr + \lambda \tilde{\vrs}) Z }{(Z+ \lambda )^2} + \frac{\lambda \tilde{\vrs}}{ Z + \lambda  }= \tilde{\vrs}$, thus \eqref{fin} becomes trivial.

{\it Case 2. } For $Z>0$,  we first notice that  $ \frac{ \vr + \lambda \tilde{\vrs}}{Z+ \lambda } \leq \max\{\tilde\vrs,\frac{1}{c_\star}\}$. By means of the strong convergence of $\vr_\lambda=\vr+\lambda$ and $Z_\lambda=Z+\lambda$ to $\vr$ and $Z$, respectively, we can now let $\lambda\to0$ in \eqref{fin} to obtain
\eqh{\pt \lr{\frac{\vr}{Z}}+\Div\lr{\frac{\vr}{Z}\vu}-\frac{\vr}{Z}\Div\vu=0.}
Obviously, $\vrs$ defined as $\frac{\vr}{Z}$ satisfies
$\vr^*\in[\min\{(c^\star)^{-1},\tilde\vrs\}, \max\{(c_\star)^{-1},\tilde\vrs\}]$ a.e. in $(0,T)\times\Omega$, and thus $Z=\frac{\vr}{\vrs}$  almost everywhere in $(0,T)\times\Omega$. This leads to the conclusion that the condition $(1-Z)\pi=0$ can be replaced by $\lr{1-\frac{\vr}{\vrs}}\pi=0$, or, equivalently by $(\vrs-\vr)\pi=0$, where the product is defined as in \eqref{product}. This finishes the proof of Theorem \ref{Th:main}. $\Box$

\section{Numerical scheme}\label{Sec:num}

Numerical simulation of two-phase flows with free boundary requires to design a method that captures phase transition and the limit behaviour. 
In our case, the main difficulty is to propose a scheme that is independent of singular pressure parameter $\ep$ \eqref{pie}. This property is referred to as the Asymptotic Preserving (AP) property, see e.g.\cite{DeJiLi2007}. The passage with $\ep\to 0$ resembles the low Mach number limit problem, where one observes incompressible behaviour in regions where the Mach number approaches $0$. However, in the contrary to the low Mach number limit, in our model the singularity is embedded in the definition of the singular pressure $\pi$.

We adapt numerical method from \cite{DeHuNa} where the Euler system with constant maximal density constraint has been studied. One dimensional version of \textit{the Direct method} is modified and extended to capture variable density constraint and the viscosity term in the momentum equation. 
In what follows we focus only on the new elements of our approach, for the detailed description of the other parts we refer to \cite{DeHuNa} and references therein. Following this rule, we present our technique on the time semi-discrete level, the discretization in space is omitted for brevity.

\subsection{Discretization scheme}
To find an approximate solution to the system \eqref{sysS} we propose a splitting algorithm. Time is discretized by  one step finite difference (with fixed time step $\Delta t$) and a finite volume method is used in space. At each time step the set of the equations is decomposed into three parts which are solved subsequently in three sub-steps.

\subsubsection{Step 1: Hyperbolic part}

The numerical solution of the Euler part of the system follows the strategy presented in  \cite{DeHuNa}. The flux in the mass balance and the singular pressure are treated implicitly:
\begin{subequations}\label{discEuler}
\begin{equation}
\frac{\vr^{n+1}-\vr^n}{\Delta t} + \Div (\vr^{n+1}\vu^{*}) = 0, \label{discEuler1}
\end{equation}
\begin{equation}
\frac{(\vr^{n+1}\vu^*)-(\vr^{n}\vu^{n}) }{\Delta t} + \Div (\vr^{n}\vu^{n} \otimes \vu^{n}) + \nabla \pi_\ep\lr{\frac{\vr^{n+1}}{\vr^{*n}}}+\nabla p\lr{\frac{\vr^{n}}{\vr^{*n}}} = \vc{0}.\label{discEuler2}
\end{equation}
\end{subequations}
The system \eqref{discEuler} is reformulated on a discrete level in terms of singular pressure. Into \eqref{discEuler1} we substitute implicit mass flux from \eqref{discEuler2} and obtain an elliptic equation for the singular pressure

\begin{equation}\label{pelip}
\vr^{n+1}(\pi_\ep) - (\Delta t)^2\Delta \pi_\ep\lr{\frac{\vr^{n+1}}{\vr^{*n}}} = \phi(\vr^n,\vr^{*n},\vu^n),
\end{equation}

where the right hand side of \eqref{pelip} reads
$$\phi(\vr^n,\vr^{*n},\vu^n) =  \rho^{n} - (\Delta t) \Div (\rho^{n}\vu^{n}) + (\Delta t)^2 \Div\left( \Div (\rho^{n}\vu^{n} \otimes \vu^{n}) + \Grad p\lr{\frac{\vr^{n}}{\vr^{*n}}}\right).$$
The singular pressure $\pi_\ep$ is computed by solving \eqref{pelip} by means of the Newton method with numerical Jacobian. In the next step we invert singular pressure to get the density. The purpose of this approach is to ensure that the density constrain is satisfied ($\vr \le \vr^*$), which now follows from the definition of $\pi_\ep$.

After the new density is obtained we directly update the momentum. This approach is called  \textit{the Direct method}, see \cite[Section 4.1]{DeHuNa}. The second approach presented in the literature is referred to as \textit{the Gauge  method} \cite{DeJiLi2007} that is based on the decomposition of the momentum $\vr\vu = \vc{a} + \Grad \varphi, \, \Div\vc{a} = 0$ into a divergence free part $\vc{a}$, and the irrotational part $\varphi$. As reported in \cite{DeHuNa} the Direct method indicates oscillations of the velocity in congested part, while the Gauge method is diffusive in uncongested region. Since the first method does not introduce any additional numerical dissipation, we adapt it for this work. For detailed description of the space discretization we refer to \cite{DeHuNa}.

We would like to emphasize that \eqref{discEuler} is a strictly hyperbolic problem, with characteristic wave speeds $\lambda_{1,2} = \vu\pm \sqrt{\frac{\partial p}{\partial (\vr/\vrs)}}$. By the definition, the Courant--Friedrichs--Lewy (CFL) condition for the
explicit part is equal to $\max(|\lambda_{1,2}|) \le \sigma \frac{\Delta x}{ \Delta t}$, with the Courant number $\sigma$. Splitting for implicit singular and explicit background pressure \eqref{discEuler2} provides that CFL condition is satisfied uniformly in $\ep$.

\subsubsection{Step 2: Diffusion }
For the sake of numerical simulations, we consider  \eqref{S} in a simplified form $\vS(\vu)= 2\mu\lap \vu$. We treat the diffusion term implicitly to avoid additional stability restrictions:
\begin{equation}\label{diffusion}
 \frac{(\vr^{n+1}\vu^{n+1})-(\vr^{n+1}\vu^*) }{\Delta t} + 2\mu\Delta\vu^{n+1} = 0.
\end{equation}
The presence of the diffusion term is important from analytical reasons only. In fact, the presented numerical scheme has been designed to solve the Euler system, therefore one can take arbitrary viscosity, such that $\mu \ge 0$. The viscosity coefficient is fixed independently of $\ep$ and small enough to recover compressible/incompressible transition. Equation \eqref{diffusion} is discretized in space by cell-centered finite volume scheme.

\subsubsection{Step 3: Congestion transport}

The transport of the congested density is  undoubtedly a main new feature of the model \eqref{sysS} and so of the presented numerical scheme.  Having the new velocity $\vu^{n+1}$ we compute the congested density as follows
$$\frac{\vr^{*n+1}-\vr^{*n}}{\Delta t} + \vu^{n+1} \nabla \vr^{*n} = 0,$$ 
where a cell-centered finite volume scheme together with upwind is used in space.

\subsection{Numerical results}

In this section we present four numerical examples that demonstrate behaviour of the proposed model in one-dimensional periodic setting. 
As a consequence of finite volume framework the proposed scheme conserves mass. 
As for domain we take the unit interval with the mesh size $\Delta x = 10^{-3}$ and the time-step $\Delta t = 10^{-4}$. In the following we choose  singular pressure parameter $\ep = 10^{-4}$ with the exponents $\alpha=\beta=2$, \eqref{pie}, and background pressure  \eqref{pbar} with the exponent $\gamma=2$, if not stated differently. 

The test cases are:

\begin{itemize}
\item Case 1 (constant congestion):
$$\begin{cases} 
\vr(x,0) &= 0.7,\\ 
\vr^*(x,0) &= 1.0,\\
\vu(x,0) &= \begin{cases}
			 0.8  \text{ if } 0.2 < x < 0.6 \\
			 -0.8 \text{ otherwise }
			\end{cases},
 \end{cases}$$
 
\begin{figure}[h]
\setlength{\unitlength}{.1\textwidth}
\begin{center}
\begin{picture}(10,10)
\put(0.5,0.5){\includegraphics[width=.95\textwidth]{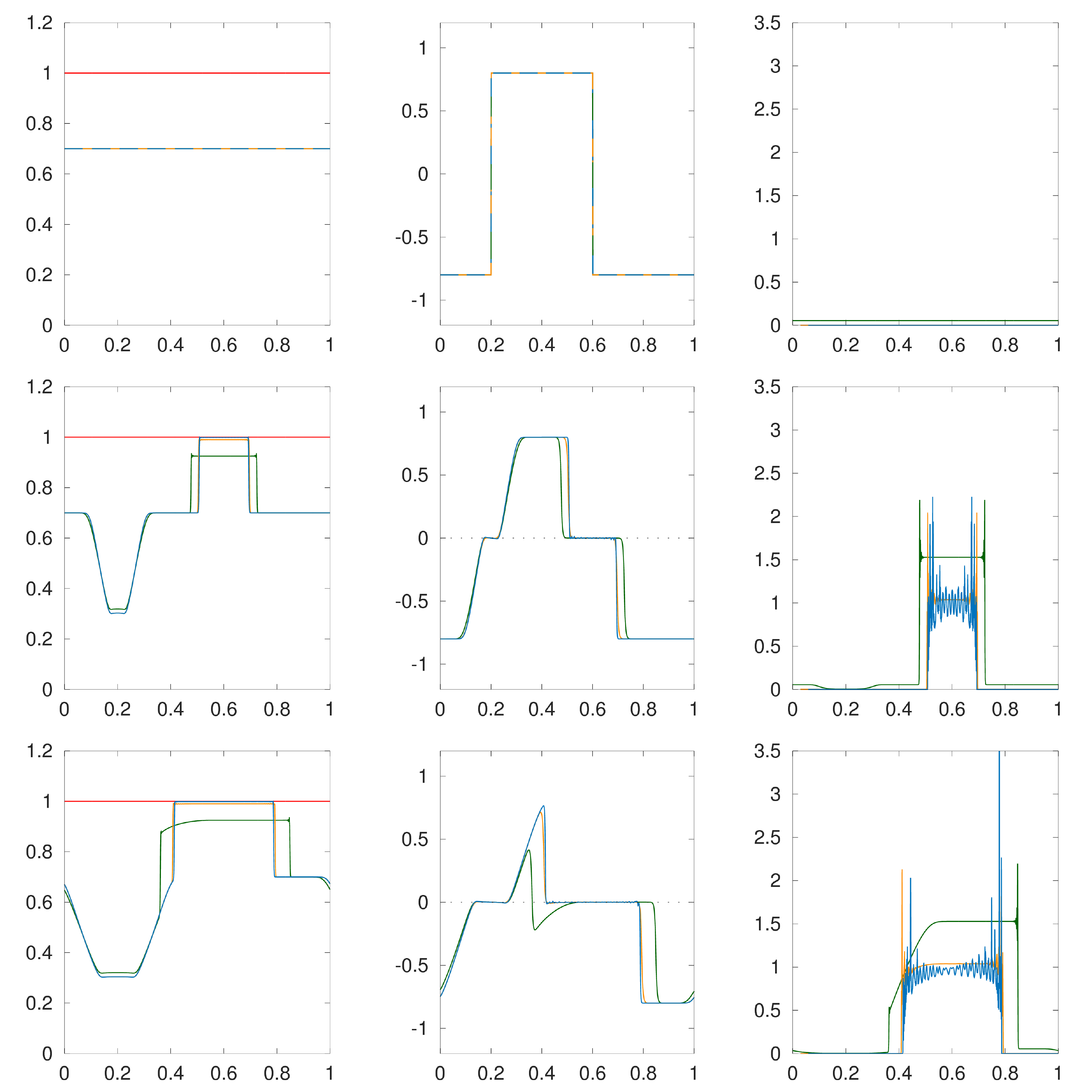}}
 \put(1.75,10.0){density}
 \put(5.0,10.0){velocity}
 \put(8.1,10.0){pressure}

 \put(0.0,1.5){\rotatebox{90}{time = 0.10}}
 \put(0.0,4.75){\rotatebox{90}{time = 0.05}}
 \put(0.0,7.85){\rotatebox{90}{time = 0.00}}

\put(2.42,7.54){\color{red}\rule{15pt}{.3pt}}
\put(2.42,4.38){\color{red}\rule{15pt}{.3pt}}
\put(2.42,1.22){\color{red}\rule{15pt}{.3pt}}

 \put(2.9,7.46){$\vr^*$}
 \put(2.9,4.3){$\vr^*$}
 \put(2.9,1.14){$\vr^*$}

\end{picture}
\end{center}
\caption{Case 1: density, velocity, and singular pressure for $\ep = 10^{-2}$(green), $\ep =10^{-4}$(yellow), $\ep =10^{-6}$(blue).}\label{fig:case1}
\end{figure}

\item Case 2:
$$ \begin{cases} 
\vr(x,0) &= 0.7,\\ 
\vr^*(x,0) &= 0.8 + 
      0.15\left(\tanh(50(x-0.4))-\tanh(50(x-0.6))\right),\\

\vu(x,0) &= \begin{cases}
			 0.8  \text{ if } 0.25 < x < 0.5 \\			 			 				-0.8  \text{ if } 0.5 < x < 0.75 \\
             0.0 \text{ otherwise }
			\end{cases},
 \end{cases} $$

\begin{figure}[h]
\setlength{\unitlength}{.1\textwidth}
\begin{center}
\begin{picture}(10,10)
\put(0.5,0.5){\includegraphics[width=.95\textwidth]{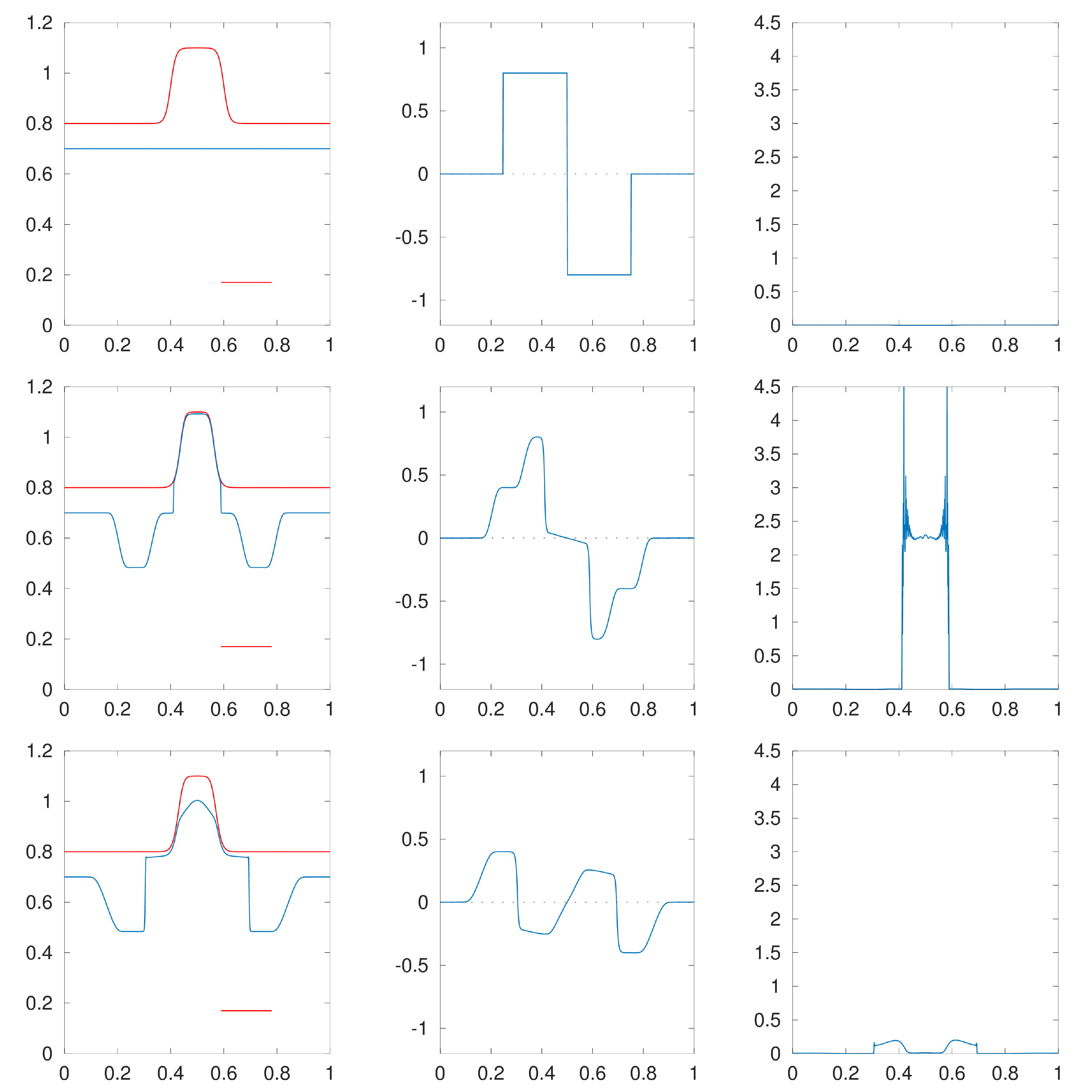}}
 \put(1.75,10.0){density}
 \put(5.0,10.0){velocity}
 \put(8.1,10.0){pressure}

 \put(0.0,1.5){\rotatebox{90}{time = 0.10}}
 \put(0.0,4.75){\rotatebox{90}{time = 0.05}}
 \put(0.0,7.85){\rotatebox{90}{time = 0.00}}

 \put(2.9,7.46){$\vr^*$}
 \put(2.9,4.3){$\vr^*$}
 \put(2.9,1.14){$\vr^*$}

\end{picture}
\end{center}
\caption{Case 2: density, velocity and singular pressure.}
\label{fig:case2}
\end{figure}

 \item Case 3:
$$\begin{cases} 
\vr(x,0) &= \begin{cases}
			 0.8  \text{ if } 0.3 < x < 0.7 \\
			 0.1 \text{ otherwise }
			\end{cases},\\ 
\vr^*(x,0) &= 0.34 + 0.3( \tanh(50(x-0.275)) - \tanh(50(x-0.725)) ),\\
\vu(x,0) &= \begin{cases}
			 0.8  \text{ if } 0.1 < x < 0.7 \\
			 0.0 \text{ otherwise }
			\end{cases},
 \end{cases}$$ 
 
\begin{figure}[h]
\setlength{\unitlength}{.1\textwidth}
\begin{center}
\begin{picture}(10,10)
\put(0.5,0.5){\includegraphics[width=.95\textwidth]{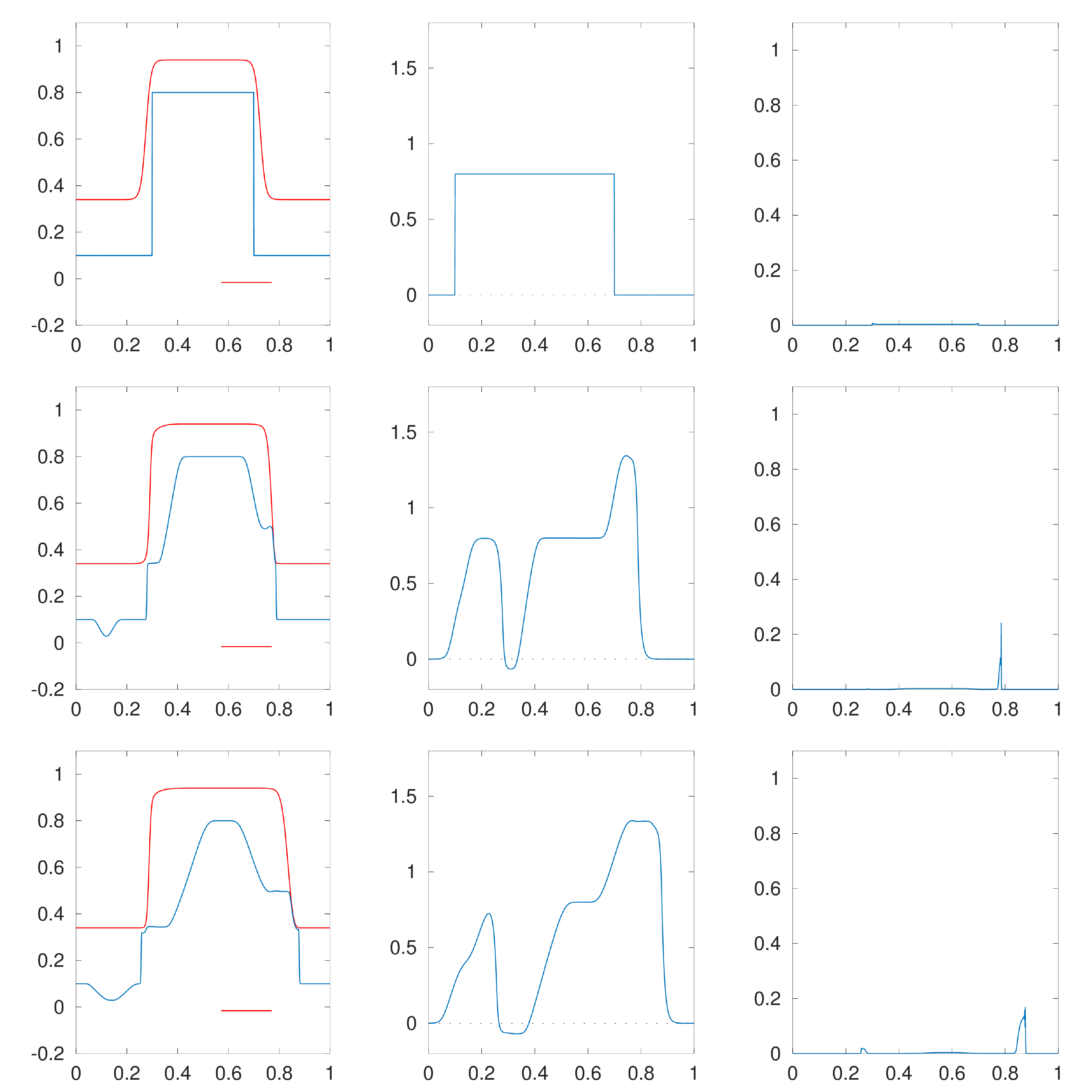}}
 \put(1.75,10.0){density}
 \put(5.0,10.0){velocity}
 \put(8.1,10.0){pressure}

 \put(0.0,1.5){\rotatebox{90}{time = 0.10}}
 \put(0.0,4.75){\rotatebox{90}{time = 0.05}}
 \put(0.0,7.85){\rotatebox{90}{time = 0.00}}

 \put(2.9,7.46){$\vr^*$}
 \put(2.9,4.3){$\vr^*$}
 \put(2.9,1.14){$\vr^*$}

\end{picture}
\end{center}
\caption{Case 3: density, velocity and singular pressure.}
\label{fig:case3}
\end{figure}
 
\item Case 4: 
 $$ \begin{cases} 
\vr(x,0) &= 0.6,\\ 
\vr^*(x,0) &= 0.9 + 0.05(\cos(10\pi x)- \cos(6\pi x) + \cos(134\pi x) + \cos(24\pi x) ),\\
\vu(x,0) &= \begin{cases}
			 0.8  \text{ if } 0.3 < x < 0.7 \\
			 -0.8 \text{ otherwise }
			\end{cases},
 \end{cases}$$ 

\begin{figure}[h]
\setlength{\unitlength}{.1\textwidth}
\begin{center}
\begin{picture}(10,13)
\put(0.5,0.5){\includegraphics[width=.95\textwidth]{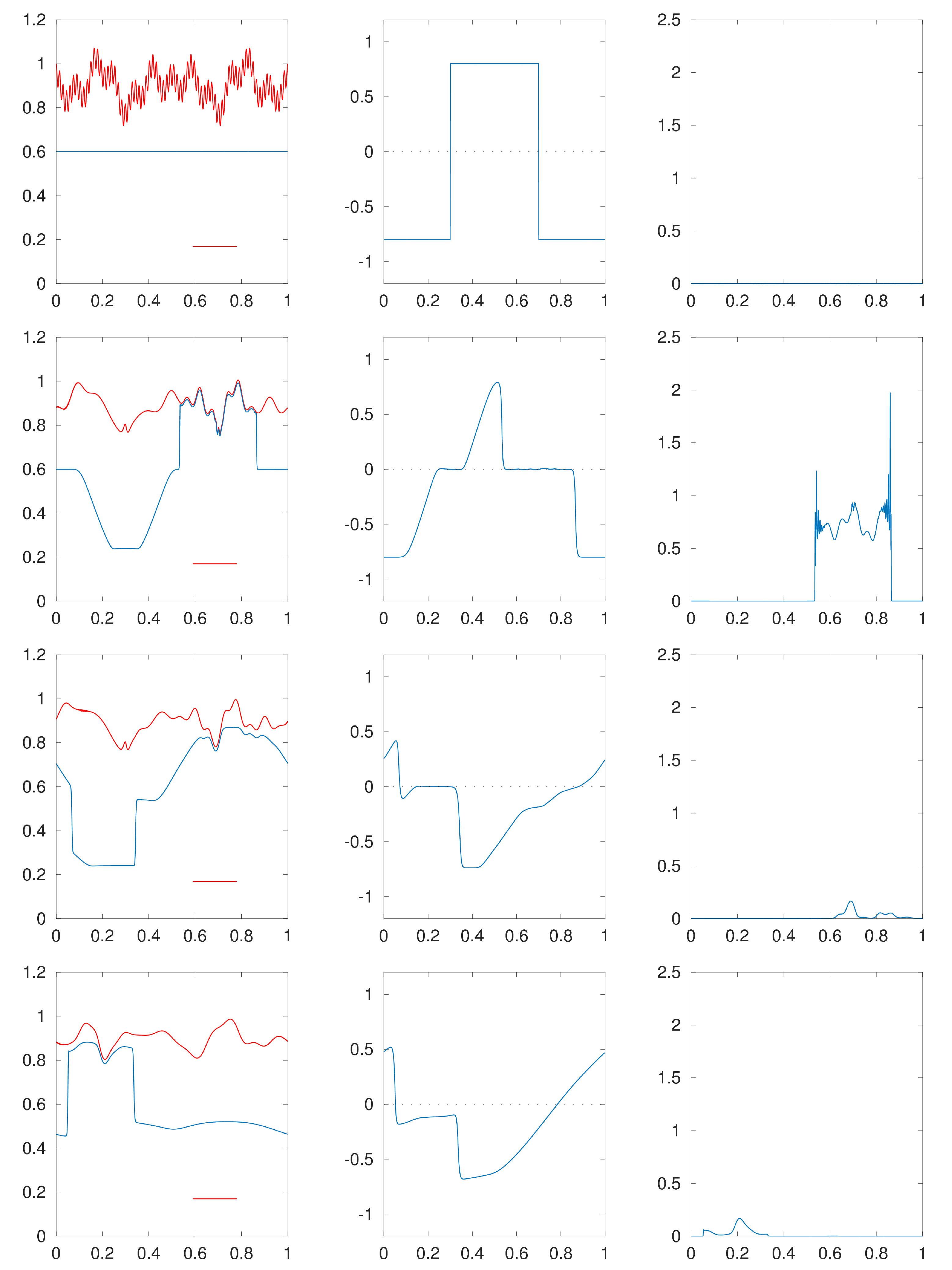}}
 \put(1.75,13.13){density}
 \put(5.0,13.13){velocity}
 \put(8.1,13.13){pressure}

 \put(0.0,1.5){\rotatebox{90}{time = 0.50}}
 \put(0.0,4.75){\rotatebox{90}{time = 0.25}}
 \put(0.0,7.85){\rotatebox{90}{time = 0.10}}
 \put(0.0,11.){\rotatebox{90}{time = 0.00}}

 \put(2.9,10.66){$\vr^*$}
 \put(2.9,7.46){$\vr^*$}
 \put(2.9,4.3){$\vr^*$}
 \put(2.9,1.14){$\vr^*$}

\end{picture}
\end{center}
\caption{Case 4: density, velocity and singular pressure.}
\label{fig:case4}
\end{figure}
\end{itemize}

Case 1 illustrates shock and rarefaction of the density for constant initial congestion density. The initial value of the congestion density stays the same for all times, due to transport. The congestions and rarefactions are created solely due to opposite initial velocities, exactly as in the analogous case from \cite{DeHuNa}.
In Figure \ref{fig:case1} we moreover  present the behaviour of unknowns for different values of parameter in the singular pressure: $\ep \in \{10^{-2},10^{-4},10^{-6}\}$. These numerical results show that the algorithm indeed satisfies the Asymptotic Preserving property. 
Note that with $\ep$ decreasing to zero we approach incompressible limit more thoroughly. However, the exact value of the maximal density constraint can never be reached by the numerically computed density.

Case 2 and Case 3 show the main feature of presented model, namely variable congestion density. For both cases the initial maximal density is set to ``smooth hat''. In the first of them the initial velocity describes the velocities of two groups of individuals that want to go in the opposite directions. The individuals close to the  contact line are more willing to compress (as the congestion density is higher). We see that initially the individuals at the rear of the groups press so intensively onto the front members, so that the whole ``hat'' is tightly filled. When this happens, we see a similar effect as for elastic collision: part of individuals start to move in the opposite direction to initially intended. 

Case 3 describes a situation, when the well organized crowd moving in one direction with the same velocity approaches a barrier ahead, being the group of individuals that move much slower and prefer to keep bigger distances between each other.
We observe how the faster individuals behind push the slower group to speed up, by filling the all available gaps between the individuals (this is where the congestion occurs at position $x\approx0.8$). This kind of behaviour could be observed, for example, at airports or in the groups of marathon runners.

Case 4 illustrates shock and rarefaction of the density when the maximal density constrain consists of a sum of cosines (periodic setting) with different frequencies. This example mimics randomness in the individual preferences of the members of population. We observe that congested regions ``freeze'' maximal density due to the zero velocity, which is consistent with the theoretical prediction.
As expected from the properties of the limiting system \eqref{pineq0}, for $\ep \ll 1$, the singular pressure $\pi_\ep$ is activated only in the congested region. In the theoretical part of the paper, this pressure is merely a nonnegative measure in the limit and our simulations seem to confirm this lack of regularity.

Another interesting feature observed in this case is the traveling wave-like behaviour of the density of the crowd.
Note that taking time derivative of \eqref{rho} and substituting $\partial_t(\rho\vu)$ from \eqref{mom} we obtain wave-like equation for density. We observe this effect on the Figure \ref{fig:case4} between time $t=0.25$ and $t=0.5$ at $x=0.2$, where two "crowds" interfere with each other. This leads to reaching the congestion density and propagation of the congestion in the opposite directions.

Looking at the Figure \ref{fig:case4} it seems that the scheme for the transport of the congestion density is quite diffusive. The high spatial frequency oscillations present at the beginning are very quickly washed away and there only subsists the small frequency components. Thus, the numerical examples presented above should be treated just as the illustration of  the behaviour of solutions to the approximation \eqref{sysS}. The thorough numerical discussion as well as the study of two-dimensional case is a purpose of our further research.

\bigskip

{\small \noindent {\bf Acknowledgements.} 
P.D. acknowledges support by the Engineering and Physical
Sciences Research Council (EPSRC) under grant no. EP/M006883/1, by the Royal
Society and the Wolfson Foundation through a Royal Society Wolfson Research
Merit Award no. WM130048 and by the National Science 
Foundation (NSF) under grant no. RNMS11-07444 (KI-Net). PD is on leave from
CNRS, Institut de Math\'ematiques de Toulouse, France. The work of E.Z. has been supported by Polish Ministry of Science and Higher Education grant "Iuventus Plus"  no. 0888/IP3/2016/74.}

{\small \noindent {\bf Data statement:} no new data were collected in the course of this research.

\end{document}